\documentclass[11pt]{article}
\input{epsf.sty}
\usepackage{hyperref}
\usepackage{amsmath, amssymb}
\newtheorem {thm}{Theorem}[section]

\def\Cox{\hfill \Box}

\def\Z{{\Bbb Z}}
\def\R{{\Bbb R}}
\def\P{{\Bbb P}}

\def\E{{\Bbb E}}

\def\0{{\bf 0}}

\def\phi{\varphi}

\def\s{\sigma}

\def\D{\Delta}

\def\L{\Lambda}

\def\T{\T}

\def\Const{\text{Const}\,}

\begin{document}

\title{A simple fluctuation lower bound for a disordered 
massless random continuous spin model in $d=2$. \thanks{Work partially
supported by   MURST (2004-06),  Cofin:  Prin 2004028108.}}
\author{
Christof
K\"ulske
\footnote{
University of Groningen, 
Department of Mathematics and Computing Sciences, 
Blauwborgje 3,   
9747 AC Groningen, 
The Netherlands
EURANDOM, LG 1.34,
\texttt{kuelske@math.rug.nl}, 
\texttt{ http://www.math.rug.nl/$\sim$kuelske/ }}\, and
Enza Orlandi
\footnote{
Dipartimento di Matematica, 
Universitˆ degli Studi "Roma Tre", 
Largo San Leonardo Murialdo, 1, 
00146 Roma,  ITALY, 
\texttt{orlandi@mat.uniroma3.it },
\texttt{http://www.mat.uniroma3.it/users/orlandi/}}
}

\maketitle

\begin{abstract}
We prove a finite volume lower bound of the order 
$\sqrt{\log N}$ on the  
delocalization of a  disordered continuous spin model  
(resp. effective interface model) in $d=2$ in a box of size $N$.  
The interaction is assumed to be massless, possibly anharmonic 
and dominated from above by a Gaussian. Disorder 
is entering via a linear source term. 
For this model delocalization with the same rate 
is proved to take place already without disorder. 
We provide a bound that is uniform in the configuration 
of the disorder, and so our proof 
shows that disorder will only enhance 
fluctuations.  \end{abstract}


\section{Introduction} \label{sect:intro}

Our model is given in terms of the formal 
infinite-volume Hamiltonian
\begin{equation}\begin{split}
&H[\eta]\left(\phi\right)
=\frac{1}{2}\sum_{i,j}p(i-j)V(\phi_i-\phi_j) - \sum_i\eta_i \phi_i\cr
\end{split}
\end{equation}
where the pair potential $V(t)$ is assumed to be twice continuously 
differentiable with an {\it upper } bound 
$V''(t)\leq c $  and  $V(t)=V(-t)$, i.e  symmetric. 
A configuration $\phi=(\phi_i)_{i \in \L}\in \R^{\L}$
can be viewed either as a continuous spin configuration 
or as an "effective interface". The $\eta=(\eta_i)_{i\in \L}$ 
denotes an arbitrary fixed configuration of external fields. 

We do not assume a lower bound on the curvature of the potential, in particular 
it might change sign and $V(t)$ might possess several minima.  
This is identical to \cite{Vel05}
and unlike to results based on the Brascamp-Lieb inequalities \cite{BL76,BrMeFr86}
which need the curvature to be strictly positive. 

Our result will be valid for all choices of the potential $V(t)$ and the random field 
configurations $\eta$ for which the integrals in finite volume are well-defined. 
For simplicity let us assume that  $V$ grows faster than linear to infinity, i.e. 
$\lim_{|x|\uparrow \infty}\frac{V(x)}{|x|^{1+\epsilon}}=\infty$. This guarantees 
the existence of the finite volume measure for all arbitrary fixed choices of $\eta\in \R^{\L}$. 

Finally $p(\cdot )$ is the transition kernel of  an aperiodic, irreducible 
random walk $X$ on $\Z^d$, assumed 
to be symmetric and, for simplicity, finite range.

Define, correspondingly the quenched finite volume Gibbs 
measures $\mu^{\hat \phi}_{N}[\eta]$, in a square $\L\equiv \L_N$ 
of sidelength $2N +1$, 
centered at the origin to be  
\begin{equation}\begin{split}\label{zwei}
&\mu^{\hat \phi}_{N}[\eta](F(\phi))\cr
&:=\frac{\int d\phi_{\L}F(\phi_{\L},\hat \phi_{\L^c})
e^{ - \frac{1}{2}\sum_{i,j\in \L}p(i-j)V(\phi_i-\phi_j) 
- \sum_{i\in\L,j\in \L^c}p(i-j)V(\phi_i-\hat\phi_j) 
+  \sum_{i\in \L}\eta_i \phi_i}}{Z_{\L}^{\hat \phi}[\eta]}\cr
\end{split}
\end{equation}
where $\hat \phi$ is a boundary condition, $\eta$ a fixed 
"frozen" configuration of random fields in $\L$ and  
$Z^{\hat \phi}_\Lambda$ is the normalization  factor.

What kind of behavior of delocalization resp. localization 
is expected to occur in a massless disordered model 
in dimension $d=2$? 
As a motivation, 
consider the Gaussian nearest neighbor 
case first, i.e. $V(x)=\frac{x^2}{2}$ 
and $p(i-j)=\frac{1}{2d}$ for $i$ and $j$ nearest neighbors. 
Then, for any fixed configuration $\eta_{\L}$, the measure 
$\mu^{\hat \phi}_{N}[\eta]$ is a Gaussian measure 
with covariance matrix $(-\D_{\L})^{-1}$ 
and expectation 
\begin{equation}\begin{split}
\int\mu^{\hat \phi}_{N}[\eta](d\phi_x)\phi_x=
 \sum_{y\in \L}(-\D_{\L})^{-1}_{x,y}\eta_y + \sum_{y\in \L^c,
|x-y|=1}(-\D_{\L})^{-1}_{x,y}\hat \phi_y.
\end{split}
\end{equation}
For every $x$ and $y$ in $\Z^d$, $d \ge 3$, the limit of $(-\D_{\L})^{-1}_{x,y}$  as $ \L \nearrow \Z^d$ exists and it is finite,  diverges like $\log N$ in $d=2$.
Taking for simplicity the random fields $\eta_y$ to be i.i.d. standard 
Gaussians, denote their expectations by $\E$,  we see that 
mean at the site $0$ of the random 
interface is itself a Gaussian variable as a linear combination 
of Gaussians and has  variance 
\begin{equation}\begin{split}
\s^2_0= \sum_{y\in \L}((-\D_{\L})^{-1}_{0,y})^2.
\end{split}
\end{equation}
This should diverge as $\int^N r (\log r)^2 d r \sim N^2(\log N)^2$ 
 when 
the sidelength $N$ of the box diverges to infinity. 
In dimension $d>2$, we have $
\int^N r^{d-1} (r^{-{(d-2)}})^2 d r $, so the interface stays 
bounded in $d>4$. 

In particular the explicit computation shows that 
delocalization is enhanced by  
randomness in the Gaussian model.    
It is however not clear whether this phenomenon 
is still present in an  anharmonic model 
where a separation  of the influence caused by the $\eta_i$'s 
is not possible and the minimizer of the Hamiltonian 
cannot be computed in a simple way. 
A priori one might not exclude the possibility that, depending on the interaction  $V$,  
a symmetrically distributed random field possibly 
stabilizes the interface.

We show in this note that this is not the case and 
the divergence is at least as strong as 
in the model without disorder, for any fixed 
field configuration. 
The method is typically two-dimensional.  Hence it does 
{\it not} show in the present form 
that  in three or four dimensions disorder will cause an 
anharmonic localized interface to diverge. 
The latter would be a continuous 
spin-analogue of the result 
in \cite{BoKu96} obtained for discrete disordered SOS-models.  
In that paper the existence of stable two-dimensional 
SOS-interfaces was excluded, using a soft martingale argument in the spirit of \cite{AW}. A disadvantage of that method however 
lies in the inability to give explicit fluctuation lower 
bounds on the behavior of the interface in finite volume. 
 
The present proof is based on a "two-dimensional" 
Mermin-Wagner type argument 
involving the entropy inequality (see \cite{Vel05}). 
The result is a quenched result, uniformly 
for all (and not only almost all) configuration of the 
disorder fields. We stress that such a "quenched instability" at any field configuration can only
hold in $d=2$, as the Gaussian interface shows. Indeed, for the Gaussian interface the instability 
of the interface is caused by fluctuations w.r.t. disorder {\it of } the groundstate, while the 
Gibbs fluctuations {\it relative to } the groundstate stay bounded. So the dimensionality of our result is correct.


\subsection{Result and proof}

\begin{thm}
Suppose $d=2$. Suppose that $\eta\in \R^{\L}$ 
is an arbitrary fixed configuration of fields. 
Then there exists a constant $c$, independent of $\eta$, such that
\begin{equation}
\label{eq:2.2general}\begin{split}
&\mu^0_N[\eta]\Bigl(|\phi_0|\geq T \sqrt{\log N}\Bigr)
\geq e^{- c T^2}.
\cr
\end{split}
\end{equation}
\end{thm}

\noindent {\bf Remark:} This generalizes the inequality of \cite{Vel05} 
to the case  of arbitrary linear disorder fields. 
We thus see that the interface is to (at least) one side "at least 
as divergent" as in the case without disorder. 

\noindent {\bf Remark 2:} Let us suppose that $\eta$ are symmetrically distributed random 
variables, possibly non-i.i.d. with {\it any } dependence structure. 
Then we have as a corollary the averaged one-sided bound 
\begin{equation}
\label{eq:2.2generalfeldmarschall}\begin{split}
&\int \P(d\eta)\mu^0_N[\eta]\Bigl( \phi_0 \geq T \sqrt{\log N}\Bigr)
\geq e^{- c T^2}/2.
\cr
\end{split}
\end{equation}
This follows immediately from the Theorem, by
symmetry. Note that no integrability assumptions on the distribution of the random fields 
are needed, given the lower bound on the potential we assume. 
\bigskip

\noindent {\bf Proof: }  
As in \cite{Vel05} we take a test-configuration $\bar\phi$, to be chosen 
later,  that 
interpolates between $\bar\phi_0=R$ and $\bar\phi_x\equiv 0$ 
for $x\in \L_N^c$.  
We define the {\it shifted  measure }Êto be  
$\mu^0_{N;\bar \phi}[\eta](\cdot)=\mu_{N}^0[\eta](\cdot + \bar\phi)$.  
Note that $\bar\phi$ does 
not depend on $\eta$. 

Let us drop the boundary 
condition from our notation and write  $\mu_N[\eta]\equiv \mu^0_N[\eta]$
in the following. 
Using the entropy-inequality we have 
\begin{equation}\begin{split}\label{imanfang}
&\mu_{N}[\eta](|\phi_0|\geq R)\cr
&=\sum_{s=\pm 1}\mu_{N}[\eta]( s \phi_0 \geq R)\cr
&=\sum_{s=\pm 1}\mu_{N}[s \eta]( \phi_0 \geq R)\cr
&=\sum_{s=\pm 1}\mu_{N;\bar \phi}[s \eta](\phi_0 \geq 0)\cr
&\geq \sum_{s=\pm 1}\mu_{N}[s\eta](\phi_0\geq 0)\exp\Bigl(
-\frac{1}{\mu_{N}[s\eta](\phi_0\geq 0)}\Bigl( 
H(\mu_{N;\bar \phi}[s \eta]|\mu_{N}[s\eta])+e^{-1}\Bigr)
\Bigr).\cr
\end{split}
\end{equation}

It remains to control the  relative entropy  
\begin{equation}\begin{split}
&H(\mu_{N;\bar \phi}[s \eta]|\mu_{N}[s\eta])
= \int\mu_{N;\bar \phi}[s \eta](d\phi)
\log\Bigl(\frac{d \mu_{N;\bar \phi}[s \eta]}{d \mu_{N}[s\eta]}
(\phi) \Bigr).
\end{split}
\end{equation}
The strategy of the proof is to 
show that we may choose $R=R(N)$ diverging with $N$ 
so that $\inf_{{\bar\phi: \bar\phi_0=R \text{ and}}\atop{\bar\phi_x\equiv 0
\text{ for } x\in \L_N^c}}
H(\mu_{N;\bar \phi}[s \eta]|\mu_{N}[s\eta])\leq \Const$, 
uniformly in $N$. 
This is identical to the case without disorder. Further we show below that 
the bound is also uniform in the field configuration $\eta$. 

Turning to the relative entropy we note 
that the appearing partition functions cancel and so 
\begin{equation}
\frac{d \mu_{N;\bar \phi}[s \eta]}{d \mu_{N}[s\eta]}
(\phi)
= \exp\Bigl( 
- H^0_\L[s\eta](\phi-\bar \phi)+ H^0_\L[s\eta](\phi)
\Bigr).
\end{equation}
Therefore
\begin{equation}\label{relent}
\begin{split}
H(\mu_{N;\bar \phi}[s \eta]|\mu_{N}[s\eta])
&=\int\mu_{N}[s\eta](d\phi) 
\Bigl(- H^0_\L[s\eta](\phi)+ H^0_\L[s\eta](\phi+\bar \phi)
\Bigr).\cr
\end{split}
\end{equation}
We rewrite the integrand of (\ref{relent}) in the form 
\begin{equation}
\begin{split}
&- H^0_\L[s\eta](\phi)+ H^0_\L[s\eta](\phi+\bar \phi)\cr
&=\frac{1}{2}\sum_{i,j\in \L}p(i-j)\Bigl(V(\phi_i-\phi_j)- 
V(\phi_i-\phi_j+\bar \phi_i-\bar \phi_j)
\Bigr)\cr
& + \sum_{i \in\L,j\in \L^c}
p(i-j)\Bigl(V(\phi_i) - V(\phi_i+\bar\phi_i)
\Bigr)  -s\sum_{i\in \L}\eta_i\bar \phi_i.
\end{split}
\end{equation}
We use now the symmetrization trick brought 
to our attention by Yvan Velenik 
(cf. \cite{Pfi81,IoShloVel02}) which here simply consists 
in estimating 
\begin{equation}\begin{split}
&H(\mu_{N;\bar \phi}[s \eta]|\mu_{N}[s\eta])
\leq \sum_{s'=\pm 1}
H(\mu_{N;\bar \phi}[s' \eta]|\mu_{N}[s'\eta]).
\end{split}
\end{equation}
We note that the $s'$-sum over the 
random potential term simply vanishes since it is independent 
of $\phi$ and hence 
\begin{equation}
\begin{split}
& \sum_{s'=\pm 1}s'\sum_{i\in \L}\eta_i\bar \phi_i=0.
\end{split}
\end{equation}
Finally, to estimate the other term we make apparent 
the quenched measure
$\frac{\mu_{N}[\eta]+ \mu_{N}[-\eta]}{2}$ and use 
its symmetry. 

So we have that 
\begin{equation}
\begin{split}
&2 \int\frac{\mu_{N}[\eta]+ \mu_{N}[-\eta]}{2}(d\phi) \Bigl(V(\phi_i-\phi_j)- 
V(\phi_i-\phi_j+\bar \phi_i-\bar \phi_j)
\Bigr)\cr
&\leq 2
\int\frac{\mu_{N}[\eta]+ \mu_{N}[-\eta]}{2}(d\phi)
V'(\phi_i-\phi_j)(\bar \phi_i-\bar \phi_j)+ c(\bar \phi_i-\bar \phi_j)^2\cr
& = c (\bar \phi_i-\bar \phi_j)^2.
\end{split}
\end{equation}
This gives 
\begin{equation}
\begin{split}
&H(\mu_{N;\bar \phi}[s \eta]|\mu_{N}[s\eta])
\leq \frac{c}{2}\sum_{i,j\in \L}p(i-j)
(\bar \phi_i-\bar \phi_j)^2+ c \sum_{i \in\L,j\in \L^c}
p(i-j)\bar \phi^2_i.
\end{split}
\end{equation}
for both $s=\pm 1$. Keeping only the $s$-term in inequality 
(\ref{imanfang}) for which 
$\mu_{N}[s\eta](\phi_0\geq 0)\geq \frac{1}{2}$ one obtains 
in fact 
\begin{equation}\begin{split}
&\mu_{N}[\eta](|\phi_0|\geq R)\cr
&\geq \frac{1}{2}\exp\Bigl(
- 2\Bigl( \frac{c}{2}\sum_{i,j\in \L}p(i-j)
(\bar \phi_i-\bar \phi_j)^2+ c \sum_{i \in\L,j\in \L^c}
p(i-j)\bar \phi^2_i +e^{-1}\Bigr)
\Bigr).\cr
\end{split}
\end{equation}
This is exactly the same bound as in the case of 
vanishing $\eta$.  
It remains to choose $\bar\phi$ optimal. Denoting by $X$ a random walk 
with the transition kernel $p$, we choose  as in   \cite{Vel05}, 
  $ \bar \phi_i = R \P_i [T_{\{0\}} < \tau_{\L_N}] $,  where $ \P_i $ is the measure 
  of the random walk started in the point $i$, 
  $T_{\{0\}}= \min \{ n: X_n=0\}$
  and   $\tau_{\L_N} =\min  \{ n: X_n \notin \L_N \} $.  
  Taking into account the estimate  \cite{Law} 
$$ \P_i[T_{\{0\}} < \tau_{\L_N}]  \simeq \frac {\ln(|i|+1)} {\ln(N+1)} $$
 gives  indeed
\begin{equation} 
\inf_{{\bar\phi: \bar\phi_0=R \text{ and}}\atop{\bar\phi_x\equiv 0
\text{ for } x\in \L_N^c}}
\Biggl(
\frac{c}{2}\sum_{i,j\in \L_N}p(i-j)
(\bar \phi_i-\bar \phi_j)^2+ c \sum_{i \in\L_N,j\in \L_N^c}
p(i-j)\bar \phi^2_i \Biggr)
\leq \Const 
\frac{R^2}{\log N}.
\end{equation}
Choosing $ R= T \sqrt{ \log N} $ one obtains \eqref {eq:2.2general}. 
$\Cox$
\bigskip 

{\bf Acknowledgements: } 
This work was stimulated by an inspiring mini-course 
of Yvan Velenik at the workshop  "Random Interfaces and Directed Polymers" in Leipzig (2005) whom we would also like 
to thank for a very useful comment on an earlier version.
CK would like to thank University Rome  Tre for hospitality.

\end{document}